\documentclass[11pt]{article}
\usepackage{amsfonts, amsmath,amssymb}
\usepackage[dvips]{graphicx}

\newcommand{\file}{$\ti{\ }$/tekst/habil/habil.tex\quad}
\renewcommand{\file}{}

\newcommand{\detail}[1]{\par\noi{\bf [Proof detail\ }{#1}
\hfill{\bf ]}\par\noi\hspace{-4pt}}
\renewcommand{\detail}[1]{}

\newcommand{\drop}[1]{}

\newcommand{\dis}{\displaystyle}
\newcommand{\txt}{\textstyle}

\newcommand{\med}{\medskip}
\newcommand{\noi}{\noindent}

\newcommand{\halmos}{\rule{1ex}{1.4ex}}
\def \qed {\nopagebreak{\hspace*{\fill}$\halmos$\medskip}}

\newtheorem{theorem}{Theorem}
\newtheorem{proposition}[theorem]{Proposition}
\newtheorem{corollary}[theorem]{Corollary}
\newtheorem{conjecture}[theorem]{Conjecture}
\newtheorem{lemma}[theorem]{Lemma}

\newtheorem{remark}[theorem]{Remark}
\newtheorem{defi}[theorem]{Definition}
\newtheorem{note}[theorem]{Remark}
\newtheorem{fact}[theorem]{Claim}
\newcommand{\bt}{\begin{theorem}}
\newcommand{\et}{\end{theorem}}
\newcommand{\bl}{\begin{lemma}}
\newcommand{\el}{\end{lemma}}
\newcommand{\bp}{\begin{proposition}}
\newcommand{\ep}{\end{proposition}}
\newcommand{\bcor}{\begin{corollary}}
\newcommand{\ecor}{\end{corollary}}
\newcommand{\br}{\begin{remark}\rm}
\newcommand{\er}{\end{remark}}
\newcommand{\bcon}{\begin{conjecture}}
\newcommand{\econ}{\end{conjecture}}
\newcommand{\bfa}{\begin{fact}}
\newcommand{\efa}{\end{fact}}
\newcommand{\bdf}{\begin{defi}\rm}
\newcommand{\edf}{\hfill$\Diamond$\end{defi}}
\newcommand{\brm}{\begin{note}\rm}
\newcommand{\erm}{\hfill$\Diamond$\end{note}}

\newcommand{\be}{\begin{equation}}
\newcommand{\ee}{\end{equation}}
\newcommand{\ba}{\begin{array}}
\newcommand{\ea}{\end{array}}
\newcommand{\bc}{\be\begin{array}{r@{\,}c@{\,}l}}
\newcommand{\ec}{\end{array}\ee}

\newcommand{\al}{\alpha}
\newcommand{\bet}{\beta}
\newcommand{\ga}{\gamma}

\newcommand{\de}{\delta}

\newcommand{\eps}{\varepsilon}

\newcommand{\La}{\Lambda}

\newcommand{\tet}{\theta}

\newcommand{\Ei}{{\cal E}}

\newcommand{\Xc}{{\cal X}}
\newcommand{\Yi}{{\cal Y}}
\newcommand{\Zi}{{\cal Z}}

\newcommand{\R}{{\mathbb R}}
\newcommand{\N}{{\mathbb N}}

\newcommand{\beh}{\backslash}

\newcommand{\Asto}[1]{\underset{{#1}\to\infty}{\Longrightarrow}}

\newcommand{\ov}{\overline}

\newcommand{\ffrac}[2]{{\textstyle\frac{{#1}}{{#2}}}}
\newcommand{\dif}[1]{\ffrac{\partial}{\partial{#1}}}
\newcommand{\diff}[1]{\ffrac{\partial^2}{{\partial{#1}}^2}}

\newcommand{\di}{\mathrm{d}}
\newcommand{\half}{{[0,\infty)}}
\newcommand{\expo}{\mbox{\large\it e}}
\newcommand{\ex}[1]{\expo^{\,\textstyle{#1}}}

\newcommand{\Pois}{{\rm Pois}}
\newcommand{\Thin}{{\rm Thin}}

\begin{document}

\makeatletter\@addtoreset{equation}{section}
\makeatother\def\theequation{\thesection.\arabic{equation}} 

\renewcommand{\labelenumi}{{\rm (\roman{enumi})}}

\title{\vspace{-3cm}\Huge Duals and thinnings of some relatives of the contact process}
\author{\Large Jan M. Swart\vspace{6pt}\\
{\' UTIA}\\
{Pod vod\'arenskou v\v e\v z\' i 4}\\
{18208 Praha 8}\\
{Czech Republic}\\
{e-mail: swart@utia.cas.cz}
}
\date{{\small\file}\large\today\\[20pt]
\parbox{400pt}{{\bf\large Abstract} This paper considers contact
processes with additional voter model dynamics. For such models,
results of Lloyd and Sudbury can be applied to find a self-duality, as
well as dualities and thinning relations with systems of random walks
with annihilation, branching, coalescence, and deaths. We show that
similar relations, which are known from the literature for certain
interacting SDE's, can be derived as local mean field limits of the
relations of Lloyd and Sudbury.}
}
\maketitle

\noi
{\it Acknowledgement.} Work sponsored by GA\v CR grant 201/06/1323.



{\setlength{\parskip}{-2pt}\tableofcontents}

\newpage

\section{Introduction}

Lloyd and Sudbury \cite{SL95,SL97,Sud00} have studied dualities for
general spin systems which have only two-spin interactions and for
which the uniform zero configuration is a trap. Based on algebraic
considerations coming from quantum theory, they associate dualities
with certain linear operators. The requirement that these operators
are the product of operators acting only locally on each site then
leads them to consider only duality functions of the form
\be\label{psieta}
\psi_\eta(x,y)=\prod_{i\in\La}\eta^{\txt x(i)y(i)},
\ee
where $\La$ is a lattice, $x,y\in\{0,1\}^\La$ are spin configurations,
and $\eta\in\R\beh\{1\}$ is a parameter. They show that there are lots
of dualities between the models they consider. Moreover, they show
that if two models are dual to the same model (albeit with a different
duality parameter), then one of these models is a thinning of the
other.

After reviewing some main results of Lloyd and Sudbury in
Section~\ref{S:LS}, in Section~\ref{S:CVP} we turn our attention to
contact processes with additional voter model dynamics. Using
Lloyd-Sudbury theory, we show that such models are self-dual, and
moreover dual to systems of random walks with annihilation, branching,
coalescence, and deaths. We also show that the latter models are
thinnings of each other, and of the contact-voter models. In
Section~\ref{S:SDE} we consider systems of interacting SDE's used in
population dynamics. More precisely, we consider a version of the
stepping stone model with selection and mutation, as well as the super
random walk with an additional quadratic killing. We show that such
systems can be derived as `local mean field limits' of contact-voter
models, and derive (mostly well-known) dualities, thinnings, and
Poissonization relations for such models as limits of the
Lloyd-Sudbury relations.

\section{Lloyd-Sudbury theory}\label{S:LS}

\subsection{Lloyd-Sudbury dualities}

Two Markov processes $X=(X_t)_{t\geq 0}$ and $Y=(Y_t)_{t\geq 0}$ with
state spaces $E_X$ and $E_Y$ are usually\footnote{In \cite{SL95,SL97},
however, the word duality is used in a much more restricted meaning.}
called {\em dual} to each other with {\em duality function}
$\psi(\cdot,\cdot)$ if
\be
E[\psi(X_t,Y_0)]=E[\psi(X_0,Y_t)]\qquad(t\geq 0),
\ee
whenever $X$ and $Y$ are independent, for arbitrary initial laws
$P[X_0\in\cdot\,]$ and $P[Y_0\in\cdot\,]$ on $E_X$ and $E_Y$,
respectively. If the functions $\{\psi(\,\cdot,y):y\in E_Y\}$ and
$\{\psi(x,\cdot\,):y\in E_X\}$ are distribution determining, such a
duality is {\em informative}.

We will in particular be interested in the case that
$E_X=E_Y=\{0,1\}^\La$, where $\La$ is a finite or countably infinite
set, and $\psi$ is of the form (\ref{psieta}). We consider Markov
processes in $\{0,1\}^\La$ which have only two-spin interactions and
for which the uniform $0$ configuration is a trap. More precisely, we
consider Markov processes $X=(X_t)_{t\geq 0}$ in $\{0,1\}^\La$ with
formal generator of the form
\be\ba{r@{}l}\label{GLS}
G_{\rm LS}f(x):=\dis\sum_{i\neq j}q(i,j)\Big\{&
\dis\ffrac{1}{2}a x(i)x(j)\{f(x-\de_i-\de_j)-f(x)\}\\[5pt]
&\dis+bx(i)(1-x(j))\{f(x+\de_j)-f(x)\}\\[5pt]
&\dis+cx(i)x(j)\{f(x-\de_j)-f(x)\}\\[5pt]
&\dis+dx(i)(1-x(j))\{f(x-\de_i)-f(x)\}\\[5pt]
&\dis+ex(i)(1-x(j))\{f(x-\de_i+\de_j)-f(x)\}\Big\},
\ec
where $a,b,c,d,e$ are nonnegative constants and $q:\La\times\La\to\R$
is a nonnegative function such that
\be\label{qdef}
q(i,j)=q(j,i)\quad\mbox{and}\quad\sum_{j:\,j\neq i}q(i,j)=1.
\ee
Here $\de_i(j):=1$ if $i=j$ and $\de_i(j):=0$ otherwise, and we adopt
the convention that sums or suprema over $i,j$ always run over $\La$,
unless stated otherwise. The letters $a,b,c,d,e$ denote annihilation,
branching, coalescence, death, and exclusion, respectively. We may
interpret $q$ as the jump rates of a continuous time Markov process on
$\La$, called the {\em underlying motion} of the interacting particle
system $X$. The processes in (\ref{GLS}) may be constructed via a
graphical representation or via their generator; see \cite{Lig85} as a
general reference. We call the Markov process $X$ with formal
generator (\ref{GLS}) the Lloyd-Sudbury model with underlying motion
kernel $q$ and parameters $a,b,c,d,e$, or shortly, the {\em
$(q,a,b,c,d,e)$-LSM}.

Since there is no spontaneous creation of particles, it is easy to see that
a $(q,a,b,c,d,e)$-LSM $X$ started in
\be
\Ei_{\rm fin}(\{0,1\}):=\big\{x\in\{0,1\}^\La:\sum_ix(i)<\infty\big\}
\ee
stays in this space. We will say that two Lloyd-Sudbury models $X$ and $X'$
are dual to each other with duality parameter $\eta\in\R\beh\{1\}$ if
\be\label{LSdu}
E\big[\prod_i\eta^{\txt X_t(i)X'_0(i)}\big]
=E\big[\prod_i\eta^{\txt X_0(i)X'_t(i)}\big].
\ee
for all deterministic initial states $X_0$ and $X'_0$. (If $X$ and
$X'$ are independent, then integrating over the initial laws we see
that (\ref{LSdu}) holds for nondeterministic initial states as well.)
Here we assume that $X_0,X'_0\in\{0,1\}^\La$ are arbitrary if $\La$ is
finite or if $|\eta|<1$. In case $\eta=-1$ and $\La$ is infinite, we
assume that either $\sum_iX_0(i)<\infty$ or $\sum_iX'_0(i)<\infty$. To
avoid convergence problems, in case $|\eta|>1$, we assume that
$|\La|<\infty$.

We cite the next proposition from \cite[formula~(9)]{Sud00}, which is
a simplification of \cite[formula~(21)]{SL95}. While the latter
applies only to symmetric models, their formula~(20) applies to
asymmetric models as well. A simplification of that formula in the
spirit of \cite[formula~(9)]{Sud00}, plus an outline of its proof, can
be found in Appendix~\ref{A:LS}.
\bp{\bf(Dualities between Lloyd-Sudbury models)}\label{P:LS}
The $(q,a,b,c,d,e)$-LSM and $(q,a',b',c',d',e')$-LSM are dual with duality
parameter $\eta\in\R\beh\{1\}$, provided that
\be\label{dupar}
a'=a+2\eta\ga,\quad b'=b+\ga,\quad c'=c-(1+\eta)\ga,\quad d'=d+\ga,
\quad e'=e-\ga,
\ee
where $\ga:=(a+c-d+\eta b)/(1-\eta)$. In particular, each model with $b>0$
and $\eta:=(d-a-c)/b\neq 1$ is self-dual with parameter $\eta$.
\ep

\subsection{Thinnings}

In this section, we recall from \cite{SL97} that if two
$(q,a,b,c,d,e)$-LSM's have the same dual (in the sense of
Proposition~\ref{P:LS}), then one is a thinning of the other. Since we
will need this further on, we define thinnings here in a somewhat
greater generality than is needed at this point.

We may interpret the space $\N^\La$ as the space of all particle
configurations $x$ on $\La$, where each site $i\in\La$ can be occupied
by $x(i)=0,1,\ldots$ particles. A {\em $u$-thinning} of a particle
configuration $x\in\N^\La$ is then obtained by independently throwing
away particles, where a particle at $i\in\La$ is kept with probability
$u(i)$. More formally, for given $x\in\N^\La$ and $u\in[0,1]^\La$, we
can choose $i_n\in\La$ such that $x=\sum_{n=1}^m\de_{i_n}$ (where $m$
is allowed to be $\infty$), and we can choose independent
$\{0,1\}$-valued random variables $\chi_n$ with
$P[\chi_n=1]=u(i_n)$. Setting
\be\label{thindef}
\Thin_u(x):=\sum_{n=1}^m\chi_n\de_{i_n}
\ee
then defines a $u$-thinning of $x$. We will usually only be interested
in the law of $\Thin_u(x)$, and use this symbol for any random
variable with the same law as in (\ref{thindef}). If $X,U$ are random
variables with values in $\N^\La$ and $[0,1]^\La$, respectively, then
we use the symbol $\Thin_U(X)$ for any random variable with the law
\be
P[\Thin_U(X)\in\cdot\,]:=\int P[U\in\di u,\ X\in\di x]P[\Thin_u(x)\in\cdot\,].
\ee
In practice, we will only be interested in the case that $U$ and $X$ are
independent. We will sometimes need the elementary relation
\be
P[\Thin_v(\Thin_u(x))\in\cdot\,]=P[\Thin_{vu}(x)\in\cdot\,]
\qquad(v,u\in[0,1]^\La,\ x\in\N^\La).
\ee
The next lemma gives suficient conditions for one Lloyd-Sudbury model
$X^1$ to be a $v$-thinning of another Lloyd-Sudbury model $X^2$.
\bl{\bf(Thinnings of Lloyd-Sudbury models)}\label{L:thin}
Consider $(q,a^k,b^k,c^k,d^k,e^k)$-LSM's $X^k$ $(k=0,1,2)$. Assume that for
$k=1,2$, $X^k$ is dual to $X^0$ with the duality function from (\ref{psieta})
and duality parameter $\eta^{k0}$. Then
\be\label{LSthin}
P[X^1_0\in\cdot\,]=P[\Thin_v(X^2_0)\in\cdot\,]
\quad\mbox{implies}\quad
P[X^1_t\in\cdot\,]=P[\Thin_v(X^2_t)\in\cdot\,]
\qquad(t\geq 0).
\ee
provided that $v:=(1-\eta^{20})/(1-\eta^{10})\in[0,1]$.
\el
{\bf Proof} This is more or less \cite[Theorem~2.1]{Sud00} (although the
formulation of thinning is a bit different there), but since the proof is
very short, we give it here. We introduce the notation
\be
\eta^{\txt x}:=\prod_i\eta^{x(i)}\qquad(\eta\in\R,\ x\in\N^\La)
\ee
Setting $\tet^{k0}:=1-\eta^{k0}$, we see that the duality between $X^k$ and
$X^0$ can be cast in the form
\be
E\big[(1-\tet^{k0})^{\txt X^0_0X^k_t}\big]
=E\big[(1-\tet^{k0})^{\txt X^0_tX^k_0}\big]\qquad(t\geq 0).
\ee
We need the relation
\be\label{thimom}
E\big[(1-\tet)^{\txt\Thin_{\tet'}(x)}\big]=E\big[(1-\tet\tet')^{\txt x}\big].
\ee
If $\tet\in[0,1]$, this is the relation
$P[\Thin_{\tet}(\Thin_{\tet'}(x))=0]=P[\Thin_{\tet\tet'}(x)=0]$. It is
not hard to show that (\ref{thimom}) holds more generally for all
$\tet\in\R$. Now if the initial laws of $X^1$ and $X^2$ are related as
in (\ref{LSthin}), then for $t\geq 0$
\be\ba{l}
\dis E[(1-\tet^{10})^{X^0_0\Thin_v(X^2_t)}]
=E[(1-\tet^{10})^{\Thin_v(X^0_0X^2_t)}]=E[(1-v\tet^{10})^{X^0_0X^2_t}]
=E[(1-\tet^{20})^{X^0_0X^2_t}]\\[5pt]
\dis=E[(1-\tet^{20})^{X^0_tX^2_0}]=E[(1-\tet^{10})^{X^0_t\Thin_v(X^2_0)}]
=E[(1-\tet^{10})^{X^0_tX^1_0}]=E[(1-\tet^{10})^{X^0_0X^1_t}].
\ec
Since this is true for each $X^0_0$, it follows that $\Thin_u(X^2_t)$ and
$X^1_t$ are equal in distribution.\qed

\noi
{\bf Remarks} I did not check the answers to the following questions:
1.\ (Perhaps not too difficult) In Lemma~\ref{L:thin}, is it
sufficient if $X^0$ is only a formal dual, i.e., is it allowed that
some of the rates $a^0,\ldots,e^0$ are negative? 2.\ (Perhaps somewhat
more tricky) Is there a general formula, in the spirit of
(\ref{dupar}), that tells us when one Lloyd-Sudbury model is a
thinning of another?

\section{Contact-voter models}\label{S:CVP}

In this section, we apply the results from the previous section to
mixtures of the contact process and the voter model. Let $q$ be as in
(\ref{qdef}) and let $r,s,m\geq 0$ be constants. By definition, the
contact-voter process with underlying motion kernel $q$ and parameters
$r,s,m$, shortly the {\em $(q,r,s,m)$-CVP}, is the Markov process
$X=(X_t)_{t\geq 0}$ in $\{0,1\}^\La$ with formal generator
\be\ba{r@{}l}\label{GCVP}
G_{\rm CVP}f(x):=\dis\sum_{i\neq j}q(i,j)\Big\{&
(r+s)x(i)(1-x(j))\{f(x+\de_j)-f(x)\}\\[-5pt]
&\dis+r(1-x(i))x(j)\{f(x-\de_j)-f(x)\}\Big\}\\[7pt]
+&\dis m\sum_ix(i)\{f(x-\de_i)-f(x)\}.
\ec
We may interpret $X_t(i)\in\{0,1\}$ as the genetic type of an organism
living at time $t\geq 0$ at the site $i\in\La$, where $1$ is the
fitter type. An organism living at site $i$ invades a site $j$ with
rate $rq(i,j)$ if it is of type $0$, and with rate $(r+s)q(i,j)$ if it
is of type $1$. In addition, organisms of type $1$ mutate with rate
$m$ to organisms of type $0$. We call $r$ the {\em resampling rate},
$s$ the {\em selection rate}, and $m$ the {\em mutation rate}. For
$r=0$, our process is a contact process with infection rate $s$ and
recovery rate $m$, while setting $s=m=0$ yields a voter model.

Contact-voter processes fall into the class of Lloyd-Sudbury models described
in the last section. Indeed, setting
\be
a=0,\quad b=r+s,\quad c=m,\quad d=r+m,\quad e=0
\ee
in the generator in (\ref{GLS}) yields the generator in
(\ref{GCVP}). We may therefore apply Proposition~\ref{P:LS} to find
duals of contact-voter models. For a given duality parameter $\eta$,
Proposition~\ref{P:LS} tells us that the Lloyd-Sudbury model with the
parameters
\be\ba{l}\label{CVdu}
a'=2\eta\ga,\quad b'=\ffrac{1}{1-\eta}s,\quad c'=m-(1+\eta)\ga,\\[5pt]
d'=m+\ffrac{\eta}{1-\eta}s,\quad e'=-\ga,
\ec
where
\be
\ga=\ffrac{\eta}{1-\eta}s-r
\ee
is dual to the $(q,r,s,m)$-CVP, provided that the rates
$a',b',c',d',e'$ are all nonnegative (and, hence, the LSM with these
rates is well-defined). The conditions $a',e'\geq 0$ are equivalent to
$\eta\in(-\infty,0]\cup\{r/(r+s)\}$, where $r/(r+s)$ is the value of
$\eta$ that yields a self-duality. The condition $b'\geq 0$ is now
trivially fulfilled, while $d'\geq 0$ is equivalent to
$-m/(s-m)\leq\eta$ if $m<s$. The condition $c'\geq 0$ leads to
somewhat messy conditions on $\eta$, but if $-1\leq\eta$, this is also
trivially fulfilled. In the latter case, we can moreover give a nice
interpretation to the duals with $\eta\leq 0$.

Let $q$ be as before and let $\eps\in[0,1]$ and $\rho,\bet,\de\geq 0$
be constants. Let $Y=(Y_t)_{t\geq 0}$ be a Markov process in
$\{0,1\}^\La$, defined by the formal generator
\be\ba{r@{}l}\label{GRW}
G_{\rm RW}f(y):=\dis\sum_{ij}q(i,j)\Big\{&\rho y(i)
\Big((1-y(j))\{f(y+\de_j-\de_i)-f(y)\}\\[-5pt]
&\dis+y(j)\{\eps f(y-\de_j-\de_i)+(1-\eps)f(y-\de_i)-f(y)\})\Big)\\[5pt]
&\dis+\bet y(i)\Big((1-y(j))\{f(y+\de_j)-f(y)\}+\eps y(j)
\{f(y-\de_j)-f(y)\})\Big)\Big\}\\[8pt]
+&\dis\de\sum_iy(i)\{f(y-\de_i)-f(y)\}.
\ec
We may describe the process $Y$ in words as follows. Particles jump
from a site $i$ to a site $j$ with rate $\rho q(i,j)$, a particle at
the site $i$ gives with rate $\bet q(i,j)$ birth to a particle at the
site $j$, and particles die with rate $\de$. If in this process, two
particles land on the same site, then they annihilate with probability
$\eps$, and coalesce with the remaining probability. We call $Y$ a
system of random walks with annihilation, branching, and coalescence,
with underlying motion kernel $q$ and parameters $\eps,\rho,\bet,\de$,
or shortly the {\em $(q,\eps,\rho,\bet,\de)$-RW}. This type of models
falls into the class of Lloyd-Sudbury models. Indeed, the
$(q,\eps,\rho,\bet,\de)$-RW is the
$(q,2\rho\eps,\bet,\rho(1-\eps)+\bet\eps+\de,\de)$-LSM.

Reinterpreting the duals in (\ref{CVdu}), we find the following result.
\bp{\bf (Duals of contact-voter processes)}\label{P:cvdu}
\med

\noi
{\bf (a)} Each $(q,r,s,m)$-CVP with $s>0$ is self-dual with duality
parameter $r/(r+s)$.\med

\noi
{\bf (b)} A $(q,r,s,m)$-CVP is dual, with duality parameter $-\eps$,
to the
$(q,\eps,r+\frac{\eps}{1+\eps}s,\frac{1}{1+\eps}s,m-\frac{\eps}{1+\eps}s)$-RW,
provided that $m\geq\frac{\eps}{1+\eps}s$. Conversely, a
$(q,\eps,\rho,\bet,\de)$-RW is dual, with duality parameter
$\eta=-\eps$, to the
$(q,\rho-\eps\bet,(1+\eps)\bet,\de+\eps\bet)$-CVP, provided that
$\rho\geq\eps\bet$.
\ep
Note that the formulas in part~(b) simplify a lot if $\eps=0$, i.e.,
if there is no annihilation.

By Lemma~\ref{L:thin}, it follows that $(q,\eps,\rho,\bet,\de)$-RW's
are thinnings of each other, and of CVP's.
\bl{\bf(Thinnings of contact-voter processes)}\label{L:cvthin}
Fix $q$ as in (\ref{qdef}), $r,m\geq 0$, $s>0$, and
$0\leq\eps\leq\eps'\leq 1$ such that
$m\geq\frac{\eps'}{1+\eps'}s$. Then\med

\noi
{\bf (a)} The
$(q,\eps,r+\frac{\eps}{1+\eps}s,\frac{1}{1+\eps}s,m-\frac{\eps}{1+\eps}s)$-RW
is a $(1+\eps)^{-1}(1+\frac{r}{s})^{-1}$-thinning of the
$(q,r,s,m)$-CVP.\med

\noi
{\bf (b)} The
$(q,\eps',r+\frac{\eps'}{1+\eps'}s,\frac{1}{1+\eps'}s,m-\frac{\eps'}{1+\eps'}s)$-RW
is a $(1+\eps)/(1+\eps')$-thinning of the
$(q,\eps,r+\frac{\eps}{1+\eps}s,\frac{1}{1+\eps}s,m-\frac{\eps}{1+\eps}s)$-RW.
\el
In particular, it follows that each
$(q,\eps,r+\frac{\eps}{1+\eps}s,\frac{1}{1+\eps}s,m-\frac{\eps}{1+\eps}s)$-RW
is a thinning of the $(q,0,r,s,m)$-RW, i.e., a system of random walks
with branching, coalescence, and deaths (with no annihilation).

\section{Local mean field limits}\label{S:SDE}

\subsection{Basic definitions}

In this section, we will argue that two types of interacting SDE's,
and one particle system, all three with applications in population
dynamics, can be obtained as `local mean field limits' of contact
voter processes and their duals. Moreover, we will investigate how the
dualities and thinning relations of Proposition~\ref{P:cvdu} and
Lemma~\ref{L:cvthin} behave under these limits. Since our main
interest is in dualities and thinnings, we do not give detailed proofs
of our limit relations, as this would get too technical, in particular
for infinite initial states. Instead, we give the main calculations,
which for finite $\La$ are almost a proof, and roughly indicate what
needs to be done to make this precise.

In order to properly define our systems and their duals, we need a few
definitions. Let $q$ be as in (\ref{qdef}). It is possible to choose
strictly positive constants $(w_i)_{i\in\La}$ in such a way that
\be
\sum_iw_i<\infty\quad\mbox{and}\quad\sum_jq(i,j)w_j\leq Kw_i\quad
\forall i\in\La
\ee
for some $K<\infty$. We fix such constants from now on and define a
{\em Liggett-Spitzer norm} and $L_1$-norm
\be
\|x\|_w:=\sum_iw_i|x(i)|,\quad|x|:=\sum_i|x(i)|\qquad(x\in\R^\La).
\ee
Our processes of interest will be Markov processes taking values in the spaces
\be
\Ei(S):=\{x\in S^\La:\|x\|_w<\infty\}\quad\mbox{and}\quad\Ei_{\rm fin}(S)
:=\{x\in S^\La:|x|<\infty\},
\ee
where $S=\N,[0,1]$, or $\half$.

First, we will consider a {\em stepping stone model} with underlying
motion kernel $q$ and parameters $r,s,m$, shortly the {\em
$(q,r,s,m)$-SSM}. This is the Markov process $\Xc$ in $\Ei([0,1])$ or
$\Ei_{\rm fin}([0,1])$ defined by the system of stochastic
differential equations (SDE)
\bc\label{SSMsde}
\di \Xc_t(i)&=&\dis\sum_jq(j,i)(\Xc_t(j)-\Xc_t(i))\,\di t
+s\Xc_t(i)(1-\Xc_t(i))\,\di t-m\Xc_t(i)\,\di t\\
&&\dis+\sqrt{2r\Xc_t(i)(1-\Xc_t(i))}\,\di B_t(i).
\ec
This model and generalizations have been considered in
\cite{SU86}. The parameters $r,s,m$ can be interpreted as resampling,
selection, and mutation rates.

Second, we will consider a {\em branching particle system} with
additional annihilation and coalescence, with underlying motion kernel
$q$ and parameters $a,b,c,d$, shortly the {\em
$(q,a,b,c,d)$-BPS}. This is the Markov process $\Yi$ in $\Ei(\N)$ or
$\Ei_{\rm fin}(\N)$ defined by the formal generator
\bc\label{SRWgen}
G_{\rm SRW}f(y)&:=&\dis\sum_{i\neq j}q(i,j)y(i)\{f(y+\de_j-\de_i)-f(y)\}\\
&&\dis+a\sum_iy(i)(y(i)-1))\{f(y-2\de_i)-f(y)\}
+b\sum_iy(i)\{f(y+\de_i)-f(y)\}\\
&&\dis+c\sum_iy(i)(y(i)-1)\{f(y-\de_i)-f(y)\}
+d\sum_iy(i)\{f(y-\de_i)-f(y)\}.
\ec
This model, with annihilation rate $a=0$, has been considered in \cite{AS05}.

Third, we will consider a {\em super random walk} with quadratic
killing, with underlying motion kernel $q$ and parameters
$\al,\bet,\ga$, where $\al,\ga>0$, $\bet\in\R$, shortly the {\em
$(q,\al,\bet,\ga)$-SRW}. This is the Markov process $\Zi$ in
$\Ei(\half)$ or $\Ei_{\rm fin}(\half)$ defined by the system of
stochastic differential equations (SDE)
\be\label{SRWsde}
\di \Zi_t(i)=\sum_jq(j,i)(\Zi_t(j)-\Zi_t(i))\,\di t+\sqrt{\al \Zi_t(i)}\,\di B_t(i)+\bet \Zi_t(i)\,\di t-\ga \Zi_t(i)^2\,\di t.
\ee
This model has been considered in \cite{HW05}.

We will approximate these processes with contact-voter models and
their duals, living on a lattice of the form
$\La^{(N)}:=\La\times\{1,\ldots,N\}$, with underlying motion kernel
\be\label{quN}
q^{(N)}\big((i,k),(j,l)\big):=\left\{\ba{ll}\frac{\nu^{(N)}}{N}q(i,j)&\mbox{if }i\neq j,\\[5pt]
\frac{1-\nu^{(N)}}{N-1}&\mbox{if }i=j,\ k\neq l.\ea\right.
\ee
We interpret the collection of sites $\{(i,k):k=1,\ldots,N\}$ as a
colony of $N$ organisms. Observe that (\ref{quN}) says that our
underlying motion jumps with rate $\nu^{(N)}$ to a site in a different
colony, and with rate $1-\nu^{(N)}$ to a different site in the same
colony. We will be interested in the case that $N\to\infty$ and
$\nu^{(N)}\to 0$, i.e., we have large colonies and most of the
interaction takes place between organisms living in the same colony.

When considering dualities, it will be useful to assume that the laws
of our approximating $\{0,1\}^{\La^{(N)}}$-valued processes $X^{(N)}$
are symmetric with respect to permutations of the individuals within a
colony, i.e.,
\be\label{colsym}
P[\Pi_{i,\pi}(X^{(N)}_t)\in\cdot\,]=P[X^{(N)}_t\in\cdot\,]\qquad(t\geq 0),
\ee
for each $i\in\La$ and for each permutation $\pi$ of $\{1,\ldots,N\}$, where
\be
\big(\Pi_{i,\pi}(x)\big)(j,k):=\left\{\ba{ll}x(i,\pi(k))&\mbox{ if }j=i,\\
x(j,k)&\mbox{ otherwise,}\ea\right.
\ee
is an operator that permutes the individuals in colony $i$. It is easy
to see that if (\ref{colsym}) holds at $t=0$, then also for all $t>0$.

\subsection{Stepping stone models}

We claim that stepping stone models and branching particle systems can
be obtained as local mean field limits of contact voter processes and
their duals.
\bfa{\bf(Stepping stone local mean field limit)}\label{C:cv}
Let $X^{(N)}$ be the $(q^{(N)},r^{(N)},s,m)$-CVP and let $Y^{(N)}$ be
the
$(q^{(N)},\eps,r^{(N)}+\frac{\eps}{1+\eps}s,\frac{1}{1+\eps}s,m
-\frac{\eps}{1+\eps}s)$-RW,
where $q^{(N)}$ is as in (\ref{quN}) and
\be\label{raco}
\nu^{(N)}:=1/(rN),\quad r^{(N)}:=rN.
\ee
for some $\eps\in[0,1]$, $r>0$, and $s,m\geq 0$. Set
\bc
\dis\ov X^{(N)}_t(i)&:=&\dis\frac{1}{N}\sum_{k=1}^NX^{(N)}_t(i,k),\\[5pt]
\dis\ov Y^{(N)}_t(i)&:=&\dis\sum_{k=1}^NY^{(N)}_t(i,k),
\ec
let $\Xc$ be the $(q,r,s,m)$-SSM, and let $\Yi$ be the
$(q,\eps r,\frac{1}{1+\eps}s,(1-\eps)r,m-\frac{\eps}{1+\eps})$-BPS. Then
\be\label{toSSM}
P[\ov X^{(N)}_0\in\cdot\,]\Asto{N}P[\Xc_0\in\cdot\,]\quad\mbox{implies}\quad
P[\ov X^{(N)}_t\in\cdot\,]\Asto{N}P[\Xc_t\in\cdot\,]\qquad(t\geq 0),
\ee
and
\be\label{toBPS}
P[\ov Y^{(N)}_0\in\cdot\,]\Asto{N}P[\Yi_0\in\cdot\,]\quad\mbox{implies}\quad
P[\ov Y^{(N)}_t\in\cdot\,]\Asto{N}P[\Yi_t\in\cdot\,]\qquad(t\geq 0).
\ee
\efa
The main ingredients for a proof of Claim~\ref{C:cv} can be found in
Appendix~\ref{A:cv}. The convergence in (\ref{toSSM}) is quite
natural, and shows that the stepping stone model is a good model for
gene frequencies of an organism living in large colonies, with small
migration between these colonies (compare the discussion in
\cite[Section~1]{SU86}).

We now investigate what happens to the dualities and thinning
relations from Proposition~\ref{P:cvdu} and Lemma~\ref{L:cvthin} in
the local mean field limit from Claim~\ref{C:cv}. For $\eps=0$, the
dualities below can be found in \cite[Theorem~1]{AS05}; for part~(b),
see also \cite[Lemma~2.1]{SU86}. The case $\eps>0$ in part~(b) was
discovered during work in progress of the present author with Siva
Athreya.
\bp{\bf(Duals of stepping stone models)}\label{P:sshdu}
Let $q$ be as in (\ref{qdef}), $r,m\geq 0$, $s>0$, and $\eps\in[0,1]$
such that $m\geq\frac{\eps}{1+\eps}s$. Then\med

\noi
{\bf (a)} The $(q,r,s,m)$-SSM is self-dual with duality function
\be\label{psiSSM}
\psi(x,x'):=\ex{-\frac{r}{s}\sum_ix(i)x'(i)}\qquad(x,x'\in\Ei([0,1])).
\ee
{\bf (b)} The $(q,r,s,m)$-SSM is dual to the
$(q,\eps r,\frac{1}{1+\eps}s,(1-\eps)r,m-\frac{\eps}{1+\eps})$-BPS,
with duality function
\be\label{psiSB}
\psi(x,y):=\prod_i\big(1-(1+\eps)x(i)\big)^{y(i)}\qquad
(x\in\Ei([0,1]),\ y\in\Ei(\N)).
\ee
where we assume that either $|x|<\infty$ or $|y|<\infty$ in case $\eps=1$.
\ep
{\bf Proof} These dualities can be verified by direct calculation, but
the point of this paper is to show that they occur as the limits of
the relations in Proposition~\ref{P:cvdu}, so we follow that road.

We start with self-duality. Let $\Xc$ and $\Xc'$ be independent
$(q,r,s,m)$-SSM's. We approximate $\Xc$ and $\Xc'$ with independent
$(q^{(N)},r^{(N)},s,m)$-CVP's $X^{(N)}$ and $X^{'(N)}$, respectively,
as in (\ref{toSSM}). By Proposition~\ref{P:cvdu}~(a),
\be\ba{l}\label{indu}
\dis E\Big[\prod_i\prod_{k=1}^N
\Big(\frac{r^{(N)}}{r^{(N)}+s}\Big)^{\txt X^{(N)}_t(i,k)X^{'(N)}_0(i,k)}\Big]\\[5pt]
\dis\qquad\qquad=E\Big[\prod_i\prod_{k=1}^N
\Big(\frac{r^{(N)}}{r^{(N)}+s}\Big)^{\txt X^{(N)}_0(i,k)X^{'(N)}_t(i,k)}\Big]
\qquad(t\geq 0).
\ec
Using (\ref{raco}), e rewrite the left-hand side of (\ref{indu}) as
\be\label{indu2}
E\Big[\prod_i\Big(1+N^{-1}\ffrac{s}{r}\Big)^{-\txt N N^{-1}
\sum_{k=1}^N X^{(N)}_t(i,k)X^{'(N)}_0(i,k)}\Big].
\ee
Let us assume that the laws of $X^{(N)}$ and $X^{'(N)}$ are symmetric with
respect to permutations of the individuals within colonies. Then, since
$X^{(N)}_t$ and $X^{'(N)}_0$ are moreover independent, it follows that
\be
N^{-1}\sum_{k=1}^N X^{(N)}_t(i,k)X^{'(N)}_0(i,k),
\ee
which is the fraction of the population in colony $i$ that belongs
both to $X^{(N)}_t$ and to $X^{'(N)}_0$, converges as $N\to\infty$ in
probability to $\Xc_t(i)\Xc'_0(i)$. Therefore, we see that the
expression in (\ref{indu2}) converges, as $N\to\infty$, to
\be
E\Big[\prod_i\ex{-\frac{s}{r}\Xc_t(i)\Xc'_0(i)}\big].
\ee
The same arguments apply to the right-hand side of (\ref{indu}), so we
find that $\Xc$ is self-dual with the duality function from
(\ref{psiSSM}).

To prove also part~(b) of the proposition, let $\Xc$ be a
$(q,r,s,m)$-SSM and let $\Yi$ be an independent $(q,\eps
r,\frac{1}{1+\eps}s,(1-\eps)r,m-\frac{\eps}{1+\eps})$-BPS. We
approximate $\Xc$ and $\Yi$ with independent
$(q^{(N)},r^{(N)},s,m)$-CVP's $X^{(N)}$ and
$(q^{(N)},\eps,r^{(N)}+\frac{\eps}{1+\eps}s,\frac{1}{1+\eps}s,m
-\frac{\eps}{1+\eps}s)$-RW's
$Y^{(N)}$, as in (\ref{toSSM}) and (\ref{toBPS}). By
Proposition~\ref{P:cvdu}~(b),
\be\ba{l}\label{indu3}
\dis E\Big[\prod_i\big(-\eps\big)^{\txt\sum_{k=1}^N
X^{(N)}_t(i,k)Y^{(N)}_0(i,k)}\Big]\\[5pt]
\dis\qquad\qquad=E\Big[\prod_i\big(-\eps\big)^{\txt\sum_{k=1}^N
X^{(N)}_t(i,k)Y^{(N)}_0(i,k)}\Big]\quad(t\geq 0).
\ec
Assuming that the laws of $X^{(N)}$ and $Y^{(N)}$ are symmetric with
respect to permutations of the individuals within a colony, we see
that $\sum_{k=1}^NX^{(N)}_t(i,k)Y^{(N)}_0(i,k)$, which is the number
of organisms in colony $i$ that belongs both to $X^{(N)}_t$ and to
$Y^{(N)}_0$, converges in probability as $N\to\infty$ to
$\Thin_{\Xc_t}(\Yi_0)(i)$. Treating the right-hand side of
(\ref{indu3}) in the same way and taking the limit $N\to\infty$, we
find that
\be\label{thidu}
E\Big[\big(1-(1+\eps))\big)^{\txt\Thin_{\Xc_t}(\Yi_0)}\Big]=E\Big[\big(1-(1+\eps))\big)^{\txt\Thin_{\Xc_0}(\Yi_t)}\Big].
\ee
By (\ref{thimom}), this means that $\Xc$ and $\Yi$ are dual with the
duality function from (\ref{psiSB}).\qed

\noi
There is also an analogon of Lemma~\ref{L:cvthin} for stepping stone
models and branching particle systems. We adopt the following
notation. For $y\in\half^\La$, we write $\Pois(y)$ to denote a Poisson
measure with intensity $y$, i.e., a $\N^\La$-valued random variable
whose components are independent Poisson distributed random variables
with mean $y(i)$. If $Y$ is random, then we use the symbol $\Pois(Y)$
for any random variable with the law
\be
P[\Pois(Y)\in\cdot\,]:=\int P[Y\in\di y]P[\Pois(y)\in\cdot\,].
\ee
\bl{\bf(Poissonizations of stepping stone models)}\label{L:ssthin}
Let $q$ be as in (\ref{qdef}), $r,m\geq 0$, $s>0$, and
$0\leq\eps\leq\eps'\leq 1$ such that $m\geq\frac{\eps'}{1+\eps'}s$. Then\med

\noi
{\bf (a)} The $(q,\eps
r,\frac{1}{1+\eps}s,(1-\eps)r,m-\frac{\eps}{1+\eps})$-BPS $\Yi$ is a
Poissonization of the $(q,r,s,m)$-SSM $\Xc$, weighted with
$(1+\eps)^{-1}\frac{r}{s}$, i.e.,
\be\label{Pois}
P[\Yi_0\in\cdot\,]=P[\Pois((1+\eps)^{-1}\ffrac{r}{s}\Xc_0)\in\cdot\,]
\quad\mbox{implies}\quad
P[\Yi_t\in\cdot\,]=P[\Pois((1+\eps)^{-1}\ffrac{r}{s}\Xc_t)\in\cdot\,].
\ee
{\bf (b)} The $(q,\eps'
r,\frac{1}{1+\eps'}s,(1-\eps')r,m-\frac{\eps'}{1+\eps'})$-BPS is a
$(1+\eps)/(1+\eps')$-thinning of the $(q,\eps
r,\frac{1}{1+\eps}s,(1-\eps)r,m-\frac{\eps}{1+\eps})$-BPS.
\el
This can be proved by taking the local mean field limits from
Claim~\ref{C:cv} in the thinning relations of
Lemma~\ref{L:cvthin}. Alternatively, these relations can be proved
directly using the dualities from Proposition~\ref{P:sshdu}. Since
this is very straightforward, we skip the proof.

\subsection{Super random walks}
 
We claim that super random walks with quadratic killing can be
obtained as local mean field limits of contact processes. (Note that
in this case, we do not need additional voter model dynamics.)
\bfa{\bf(Super-RW local mean field limit)}\label{C:srw}
Let $\al,\ga>0$, $\bet\in\R$, and let $X^{(N)}$ be the
$(q^{(N)},0,s^{(N)},m^{(N)})$-CVP, where $q^{(N)}$ is as in (\ref{quN}), and
\be\label{raco2}
s^{(N)}:=\sqrt{\al\ga N},\quad\nu^{(N)}:=1/s^{(N)},\quad m^{(N)}:=s^{(N)}-\bet.
\ee
Set $\phi:=\sqrt{\al/\ga}$ and 
\be
\ov X^{(N)}_t(i):=\frac{\phi}{\sqrt N}\sum_{k=1}^NX^{(N)}_t(i,k).
\ee
Let $\Zi$ be the $(q,\al,\bet,\ga)$-SRW. Then
\be\label{XtoSRW}
P[\ov X^{(N)}_0\in\cdot\,]\Asto{N}P[\Zi^0_0\in\cdot\,]\quad\mbox{implies}\quad
P[\ov X^{(N)}_t\in\cdot\,]\Asto{N}P[\Zi^0_t\in\cdot\,]\qquad(t\geq 0).
\ee
\efa
The main idea of the proof of Claim~\ref{C:srw} can be found in
Appendix~\ref{A:srw}.

The next proposition shows what happens to the well-known self-duality
of the contact process in the local mean field limit from Claim~\ref{C:srw}.
\bp{\bf(Self-duality of SRW)}\label{P:SRWdu}
The $(q,\al,\bet,\ga)$-SRW is self-dual with duality function
\be
\psi(z,z'):=\ex{-\frac{\ga}{\al}\sum_iz(i)z'(i)}.
\ee
\ep
{\bf Proof} We will show that this duality is the local mean field
limit of the contact process self-duality. Let $\Zi$ and $\Zi'$ be
independent $(q,\al,\bet,\ga)$-SRW's. Approximate $\Zi$ and $\Zi'$
with independent $(q^{(N)},0,s^{(N)},m^{(N)})$-CVP's $X^{(N)}$ and
$X^{'(N)}$ as in Claim~\ref{C:srw}. By Proposition~\ref{P:cvdu}~(a),
$X^{(N)}$ and $X^{'(N)}$ are dual with duality parameter $0$, i.e.,
with duality function
\be
\psi(x,x')=\prod_i\prod_{k=1}^N 0^{\txt x(i,k)x'(i,k)}
=\prod_i1_{\{\sum_{k=1}^Nx(i,k)x'(i,k)=0\}},
\ee
i.e.,
\be\ba{l}\label{aprdu}
\dis E\Big[\prod_i1_{\txt\{\sum_{k=1}^N X^{(N)}_t(i,k)
X^{'(N)}_0(i,k)=0\}}\Big]\\[10pt]
\dis\qquad\qquad=E\Big[\prod_i1_{\txt\{\sum_{k=1}^N
X^{(N)}_0(i,k)X^{'(N)}_t(i,k)=0\}}\Big]\qquad(t\geq 0).
\ec
Assuming that the laws of $X^{(N)}$ and $X^{'(N)}$ are symmetric with
respect to permutations of the individuals within colonies, we observe that
\be
\sum_{k=1}^N X^{(N)}_0(i,k)X^{'(N)}_t(i,k)
\ee
converges in probability to $\Pois(\frac{\ga}{\al}\Zi_0\Zi'_t)(i)$,
since for large $N$, the first process contains approximately
$\sqrt{\frac{\ga}{\al}}\Zi_0(i)N^{1/2}$ individuals, each having an
approximate probability $\sqrt{\frac{\ga}{\al}}\Zi'_t(i)N^{-1/2}$ to
belong to the dual process as well. It follows that the left-hand side
of (\ref{aprdu}) converges, as $N\to\infty$, to
\be
P[\Pois(\ffrac{\ga}{\al}\Zi_0\Zi'_t)=0]
=E\big[\ex{-\frac{\ga}{\al}\sum_i\Zi_0(i)\Zi'_t(i)}\big].
\ee
Treating the right-hand side of (\ref{aprdu}) in the same way, we find that
$\Zi$ is self-dual with the duality function from (\ref{psiSSM}).\qed

\noi
We have just seen that Proposition~\ref{P:SRWdu} follows by taking the
local mean field limit of Proposition~\ref{P:cvdu}~(a). One may wonder
if taking the same limit in Proposition~\ref{P:cvdu}~(b), one may
discover other duals of super random walks. It turns out that this is
not the case. Rather, the dual processes from
Proposition~\ref{P:cvdu}~(b) also converge to super random walks, and
one finds duality relations that, using scaling, can all be reduced to
the self-duality in Proposition~\ref{P:SRWdu}. Likewise, the thinning
relations from Lemma~\ref{L:cvthin} converge, under the local mean
field limit from Claim~\ref{C:srw}, to trivial scaling
relations. Since this is quite lengthy and does not yield anything
new, we omit the details.

\appendix

\section{Lloyd-Sudbury dualities}\label{A:LS}

\newcommand{\oldga}{\gamma}

In this appendix, we outline the proof of a generalization of
Proposition~\ref{P:LS}. Generalizing (\ref{GLS}), we look at models
with generators of the form 
\be\ba{r@{}l}\label{GLSg}
G_{\rm LS}f(x):=\dis\sum_{i\neq j}\Big\{&
\dis\ffrac{1}{2}a(i,j)x(i)x(j)\{f(x-\de_i-\de_j)-f(x)\}\\[-3pt]
&\dis+b(i,j)x(i)(1-x(j))\{f(x+\de_j)-f(x)\}\\[5pt]
&\dis+c(i,j)x(i)x(j)\{f(x-\de_j)-f(x)\}\\[5pt]
&\dis+d(i,j)x(i)(1-x(j))\{f(x-\de_i)-f(x)\}\\[5pt]
&\dis+e(i,j)x(i)(1-x(j))\{f(x-\de_i+\de_j)-f(x)\}\Big\},
\ec
where $a,b,c,d,e:\La\times\La\to\R$ are nonnegative functions such that
\be\label{supsum}
\sup_i\sum_jg(i,j)<\infty,\quad\sup_j\sum_ig(i,j)<\infty\qquad(g=a,b,c,d,e).
\ee
Without loss of generality, we may assume (as we do from now on) that
$a(i,j)=a(j,i)$. Generalizing our earlier definition, within this
appendix, let us call the Markov process with generator $G_{LS}$ the
{\em $(a,b,c,d,e)$-LSM}. We will outline a proof of the following
generalization of Proposition~\ref{P:LS}.
\bp{\bf(Dualities between Lloyd-Sudbury models)}\label{P:LS2}
The $(a,b,c,d,e)$-LSM and $(a',b',c',d',e')$-LSM are dual with duality
parameter $\eta\in\R\beh\{1\}$ provided that
\bc\label{newrates}
a'(i,j)&=&\dis a(i,j)+\eta\big(\ga(i,j)+\ga(j,i)\big),\\
b'(i,j)&=&\dis b(j,i)+\ga(i,j),\\
c'(i,j)&=&\dis c(i,j)-\ga(j,i)-\eta\ga(i,j)+\big(e(i,j)-e(j,i)\big)
+\eta\big(b(i,j)-b(j,i)\big),\\
d'(i,j)&=&\dis d(i,j)+\ga(i,j)+\big(e(i,j)-e(j,i)\big),\\
e'(i,j)&=&\dis e(j,i)-\ga(i,j),
\ec
where
\be\label{lad}
\ga(i,j):=\big(a(i,j)+c(j,i)-d(i,j)+\eta b(j,i)-e(i,j)+e(j,i)\big)/(1-\eta)
\ee
\ep
Observe that these formulas simplify to (\ref{dupar}) if the functions
$a,\ldots,e$ are symmetric. We note that an equivalent set of equations for
(\ref{newrates}) is
\bc\label{eqeq}
a'(i,j)+\eta\big(e'(i,j)+e'(j,i)\big)&=&
\dis a(i,j)+\eta\big(e(i,j)+e(j,i)\big),\\
b'(i,j)+e'(i,j)&=&\dis b(j,i)+e(j,i),\\
\ga'(i,j)&=&-\ga(j,i),\\
d'(i,j)+e'(i,j)&=&\dis d(i,j)+e(i,j),\\
e'(i,j)-\ffrac{1}{2}\ga'(j,i)&=&\dis e(j,i)-\ffrac{1}{2}\ga(i,j),
\ec
where $\ga'(i,j)$ is defined as in (\ref{lad}), with the rates $a,\ldots,e$
replaced by $a',\ldots,e'$.\med

\noi
{\bf Proof of Proposition~\ref{P:LS2}} By standard theory
\cite[Section~4.4]{EK}, we need to check that
\be\label{duco}
G_{\rm LS}\psi_\eta(\,\cdot,y)(x)=G'_{\rm LS}\psi_\eta(x,\cdot\,)(y)
\qquad(x,y\in\{0,1\}^\La),
\ee
where $G'_{\rm LS}$ is defined as in (\ref{GLSg}), with the rates
$a,\ldots,e$ replaced by $a',\ldots,e'$. Diligent calculation reveals that
\detail{
\be\ba{r@{}l}
\dis G_{\rm LS}\psi_\eta(\,\cdot,y)(x)=\sum_{i\neq j}\Big\{&
\dis\ffrac{1}{2}a(i,j)x(i)x(j)(\eta^{-y(i)-y(j)}-1)\\[5pt]
&\dis+b(i,j)x(i)(1-x(j))(\eta^{y(j)}-1)\\[5pt]
&\dis+c(i,j)x(i)x(j)(\eta^{-y(j)}-1)\\[5pt]
&\dis+d(i,j)x(i)(1-x(j))(\eta^{-y(i)}-1)\\[5pt]
&\dis+e(i,j)x(i)(1-x(j))(\eta^{-y(i)+y(j)}-1)\Big\}\psi_\eta(x,y).
\ec
The terms in brackets can be rewritten as
\be\ba{l}
\ffrac{1}{2}a(i,j)x(i)x(j)\Big((\eta^{-2}-1)y(i)y(j)+(\eta^{-1}-1)
\big(y(i)(1-y(j))+(1-y(i))y(j)\big)\Big)\\[5pt]
+b(i,j)x(i)(1-x(j))(\eta-1)y(j)\\[5pt]
+c(i,j)x(i)x(j)(\eta^{-1}-1)y(j)\\[5pt]
+d(i,j)x(i)(1-x(j))(\eta^{-1}-1)y(i)\\[5pt]
+e(i,j)x(i)(1-x(j))\big((\eta^{-1}-1)y(i)(1-y(j))+(\eta-1)(1-y(i))y(j)\big),
\ec
yielding
\be\ba{l}
\dis G_{\rm LS}\psi_\eta(\,\cdot,y)(x)\\[5pt]
\dis=\sum_{i\neq j}\Big\{\Big(\ffrac{1}{2}a(i,j)
\big((\eta^{-2}-1)-2(\eta^{-1}-1)\big)
+e(i,j)\big((\eta^{-1}-1)+(\eta-1)\big)\Big)x(i)x(j)y(i)y(j)\\[5pt]
\dis\quad+\Big(\ffrac{1}{2}a(i,j)(\eta^{-1}-1)-d(i,j)(\eta^{-1}-1)-e(i,j)
(\eta^{-1}-1)\Big)x(i)x(j)y(i)\\[5pt]
\dis\quad+\Big(\ffrac{1}{2}a(i,j)(\eta^{-1}-1)-b(i,j)(\eta-1)+c(i,j)
(\eta^{-1}-1)-e(i,j)(\eta-1)\Big)x(i)x(j)y(j)\\[5pt]
\dis\quad+\Big(-e(i,j)\big((\eta^{-1}-1)+(\eta-1)\big)\Big)x(i)y(i)y(j)\\[5pt]
\dis\quad+\Big(d(i,j)(\eta^{-1}-1)+e(i,j)(\eta^{-1}-1)\Big)x(i)y(i)\\[5pt]
\dis\quad+\Big(b(i,j)(\eta-1)+e(i,j)(\eta-1)\Big)x(i)y(j)\Big\}\psi_\eta(x,y).
\ec
Since $(\eta^{-2}-1)-2(\eta^{-1}-1)=(\eta^{-1}-1)(\eta^{-1}+1)
-2(\eta^{-1}-1)=(\eta^{-1}-1)\big((\eta^{-1}+1)-2\big)=(\eta^{-1}-1)^2$
and $(\eta^{-1}-1)+(\eta-1)=(1-\eta)(\eta^{-1}-1)$ this can be simplified to
\be\ba{l}
\dis G_{\rm LS}\psi_\eta(\,\cdot,y)(x)\\[5pt]
\dis=(\eta^{-1}-1)\sum_{i\neq j}\Big\{\Big(\ffrac{1}{2}(\eta^{-1}-1)a(i,j)
+(1-\eta)e(i,j)\Big)x(i)x(j)y(i)y(j)\\[5pt]
\dis\quad+\Big(\ffrac{1}{2}a(i,j)-d(i,j)-e(i,j)\Big)x(i)x(j)y(i)\\[5pt]
\dis\quad+\Big(\ffrac{1}{2}a(i,j)+\eta b(i,j)+c(i,j)
+\eta e(i,j)\Big)x(i)x(j)y(j)\\[5pt]
\dis\quad-(1-\eta)e(i,j)x(i)y(i)y(j)\\[5pt]
\dis\quad+\Big(d(i,j)+e(i,j)\Big)x(i)y(i)\\[5pt]
\dis\quad-\eta\Big(b(i,j)+e(i,j)\Big)x(i)y(j)\Big\}\psi_\eta(x,y).
\ec
Introducing some arbitrary ordering on $\La$, we have terms with $i<j$ and
$i>j$. Relabeling the indices in the terms with $i>j$, using the fact that
$a(i,j)=a(j,i)$, we obtain
\be\ba{l}
\dis G_{\rm LS}\psi_\eta(\,\cdot,y)(x)\\[5pt]
\dis=(\eta^{-1}-1)\sum_{i<j}\Big\{\Big((\eta^{-1}-1)a(i,j)+(1-\eta)
\big(e(i,j)+e(j,i)\big)\Big)x(i)x(j)y(i)y(j)\\[5pt]
\dis\quad+\Big(a(i,j)-d(i,j)-e(i,j)+\eta b(j,i)+c(j,i)+\eta e(j,i)\Big)
x(i)x(j)y(i)\\[5pt]
\dis\quad+\Big(a(j,i)-d(j,i)-e(j,i)+\eta b(i,j)+c(i,j)+\eta e(i,j)\Big)
x(i)x(j)y(j)\\[5pt]
\dis\quad-(1-\eta)e(i,j)x(i)y(i)y(j)\\[5pt]
\dis\quad-(1-\eta)e(j,i)x(j)y(i)y(j)\\[5pt]
\dis\quad+\Big(d(i,j)+e(i,j)\Big)x(i)y(i)\\[5pt]
\dis\quad+\Big(d(j,i)+e(j,i)\Big)x(j)y(j)\\[5pt]
\dis\quad-\eta\Big(b(i,j)+e(i,j)\Big)x(i)y(j)\\[5pt]
\dis\quad-\eta\Big(b(j,i)+e(j,i)\Big)x(j)y(i)\Big\}\psi_\eta(x,y),
\ec
which can be rewritten as}
\be\ba{l}\label{lidu}
\dis G_{\rm LS}\psi_\eta(\,\cdot,y)(x)\\[5pt]
\dis=(\eta^{-1}-1)\sum_{i\neq j}\Big\{\ffrac{1}{2}\Big((\eta^{-1}-1)a(i,j)
+(1-\eta)\big(e(i,j)+e(j,i)\big)\Big)x(i)x(j)y(i)y(j)\\[5pt]
\dis\quad+\Big(a(i,j)-d(i,j)-e(i,j)+\eta b(j,i)+c(j,i)+\eta e(j,i)\Big)
x(i)x(j)y(i)\\[5pt]
\dis\quad-(1-\eta)e(i,j)x(i)y(i)y(j)\\[5pt]
\dis\quad+\Big(d(i,j)+e(i,j)\Big)x(i)y(i)\\[5pt]
\dis\quad-\eta\Big(b(i,j)+e(i,j)\Big)x(i)y(j)\Big\}\psi_\eta(x,y).
\ec
and
\detail{Interchanging the roles of $x$ and $y$ and replacing $a,\ldots,e$
by $a',\ldots,e'$, we see that
\be\ba{l}
\dis G'_{\rm LS}\psi_\eta(x,\cdot\,)(y)\\[5pt]
\dis=(\eta^{-1}-1)\sum_{i\neq j}\Big\{\ffrac{1}{2}\Big((\eta^{-1}-1)a'(i,j)
+(1-\eta)\big(e'(i,j)+e'(j,i)\big)\Big)x(i)x(j)y(i)y(j)\\[5pt]
\dis\quad+\Big(a'(i,j)-d'(i,j)-e'(i,j)+\eta b'(j,i)+c'(j,i)+\eta e'(j,i)\Big)
x(i)y(i)y(j)\\[5pt]
\dis\quad-(1-\eta)e'(i,j)x(i)x(j)y(i)\\[5pt]
\dis\quad+\Big(d'(i,j)+e'(i,j)\Big)x(i)y(i)\\[5pt]
\dis\quad-\eta\Big(b'(i,j)+e'(i,j)\Big)x(j)y(i)\Big\}\psi_\eta(x,y).
\ec
In order to make the terms match up with (\ref{lidu}), we must relabel
indices in the last term and change the order of the terms. The result is}
\be\ba{l}\label{redu}
\dis G'_{\rm LS}\psi_\eta(x,\cdot\,)(y)\\[5pt]
\dis=(\eta^{-1}-1)\sum_{i\neq j}\Big\{\ffrac{1}{2}\Big((\eta^{-1}-1)a'(i,j)
+(1-\eta)\big(e'(i,j)+e'(j,i)\big)\Big)x(i)x(j)y(i)y(j)\\[5pt]
\dis\quad-(1-\eta)e'(i,j)x(i)x(j)y(i)\\[5pt]
\dis\quad+\Big(a'(i,j)-d'(i,j)-e'(i,j)+\eta b'(j,i)+c'(j,i)+\eta e'(j,i)\Big)
x(i)y(i)y(j)\\[5pt]
\dis\quad+\Big(d'(i,j)+e'(i,j)\Big)x(i)y(i)\\[5pt]
\dis\quad-\eta\Big(b'(j,i)+e'(j,i)\Big)x(i)y(j)\Big\}\psi_\eta(x,y).
\ec
Comparing (\ref{lidu}) and (\ref{redu}), we see that (\ref{duco}) is
satisfied provided that
\be\ba{l}\label{ducofin}
\dis a'(i,j)+\eta\big(e'(i,j)+e'(j,i)\big)=a(i,j)
+\eta\big(e(i,j)+e(j,i)\big),\\[5pt]
\dis-(1-\eta)e'(i,j)=a(i,j)-d(i,j)-e(i,j)+\eta b(j,i)+c(j,i)
+\eta e(j,i),\\[5pt]
\dis a'(i,j)-d'(i,j)-e'(i,j)+\eta b'(j,i)+c'(j,i)+\eta e'(j,i)
=-(1-\eta)e(i,j),\\[5pt]
\dis d'(i,j)+e'(i,j)=d(i,j)+e(i,j),\\[5pt]
\dis b'(j,i)+e'(j,i)=b(i,j)+e(i,j).
\ec
which can be solved to yield (\ref{newrates}).\qed

\detail{Indeed, setting
\bc
\oldga(i,j)&:=&\big(a(i,j)+c(j,i)-d(i,j)+\eta b(j,i)\big)/(1-\eta),\\
\zeta(i,j)&:=&\big(e(i,j)-e(j,i)\big)/(1-\eta),
\ec
and using the fact that $e(j,i)+\zeta(i,j)=\big((1-\eta e(i,j)
+\eta(e(i,j)-e(j,i))\big)/(1-\eta)=\big(e(i,j)-\eta e(j,i)\big)/(1-\eta)
=\big((1-\eta)e(j,i)+(e(i,j)-e(j,i))\big)=e(j,i)+\zeta(i,j)$, the equations
(\ref{ducofin}) may be rewritten as 
\be\ba{l}\label{ducofin2}
\dis a'(i,j)=a(i,j)+\eta\big(e(i,j)-e'(i,j)\big)
+\eta\big(e(j,i)-e'(j,i)\big),\\[5pt]
\dis e'(i,j)=e(i,j)+\eta\zeta(i,j)-\oldga(i,j)\\[5pt]
\dis\phantom{e'(i,j)}=e(j,i)+\zeta(i,j)-\oldga(i,j),\\[5pt]
\dis c'(j,i)=-a'(i,j)+d'(i,j)-\eta b'(j,i)-\big(e(i,j)-e'(i,j)\big)
+\eta(e(i,j)-e'(j,i)),\\[5pt]
\dis d'(i,j)=d(i,j)+e(i,j)-e'(i,j),\\[5pt]
\dis b'(j,i)=b(i,j)+e(i,j)-e'(j,i).
\ec
The second equation yields
\bc\label{edif}
e(i,j)-e'(i,j)&=&\oldga(i,j)-\eta\zeta(i,j),\\[5pt]
e(j,i)-e'(i,j)&=&\oldga(i,j)-\zeta(i,j).
\ec
Inserting this into the first, fourth, and fifth equation yields
\be\ba{l}\label{ducoabd}
\dis a'(i,j)=a(i,j)+\eta\big(\oldga(i,j)+\oldga(j,i)\big),\\[5pt]
\dis d'(i,j)=d(i,j)+\oldga(i,j)-\eta\zeta(i,j),\\[5pt]
\dis b'(j,i)=b(i,j)+\oldga(j,i)-\zeta(j,i).
\ec
Finally, the third equation now yields
\bc
c'(j,i)&=&\dis -a(i,j)-\eta\big(\oldga(i,j)+\oldga(j,i)\big)\\[5pt]
&&\dis +d(i,j)+\oldga(i,j)-\eta\zeta(i,j)\\[5pt]
&&\dis -\eta\big(b(i,j)+\oldga(j,i)-\zeta(j,i)\big)\\[5pt]
&&\dis -\oldga(i,j)+\eta\zeta(i,j)+\eta\big(\oldga(j,i)-\zeta(j,i)\big)\\[5pt]
&=&\dis c(j,i)-\big(a(i,j)+c(j,i)-d(i,j)+\eta b(j,i)\big)
-\eta\big(b(i,j)-b(j,i)\big)\\[5pt]
&&\dis -\eta\big(\oldga(i,j)+\oldga(j,i)\big)\\[5pt]
&=&\dis c(j,i)-(1-\eta)\oldga(i,j)-\eta\big(\oldga(i,j)
+\oldga(j,i)\big)-\eta\big(b(i,j)-b(j,i)\big).
\ec
where we have used that $\zeta(j,i)=-\zeta(i,j)$. Collecting our
results, we have found that
\bc
\dis a'(i,j)&=&\dis a(i,j)+\eta\big(\oldga(i,j)+\oldga(j,i)\big),\\[5pt]
\dis b'(j,i)&=&\dis b(i,j)+\oldga(j,i)-\zeta(j,i),\\[5pt]
\dis c'(j,i)&=&\dis c(j,i)-\oldga(i,j)-\eta\oldga(j,i)
-\eta\big(b(i,j)-b(j,i)\big),\\[5pt]
\dis d'(i,j)&=&\dis d(i,j)+\oldga(i,j)-\eta\zeta(i,j),\\[5pt]
\dis e'(i,j)&=&\dis e(j,i)+\zeta(i,j)-\oldga(i,j),
\ec
which may be rewritten as
\bc\label{acrat}
a'(i,j)&=&\dis a(i,j)+\eta\big(\oldga(i,j)+\oldga(j,i)\big),\\
b'(i,j)&=&\dis b(j,i)+\oldga(i,j)-\zeta(i,j),\\
c'(i,j)&=&\dis c(i,j)-\oldga(j,i)-\eta\oldga(i,j)
+\eta\big(b(i,j)-b(j,i)\big),\\
d'(i,j)&=&\dis d(i,j)+\oldga(i,j)-\eta\zeta(i,j),\\
e'(i,j)&=&\dis e(j,i)-\oldga(i,j)+\zeta(i,j),\\
\ec
where
\bc\label{gazet}
\oldga(i,j)&=&\big(a(i,j)+c(j,i)-d(i,j)+\eta b(j,i)\big)/(1-\eta),\\
\zeta(i,j)&=&\big(e(i,j)-e(j,i)\big)/(1-\eta).
\ec
In an attempt to further simplify this result, we put
\be
\ga(i,j):=\oldga(i,j)-\zeta(i,j).
\ee
Now formula (\ref{acrat}) can be cast in the form
\bc\label{acrat2}
a'(i,j)&=&\dis a(i,j)+\eta\big(\ga(i,j)+\ga(j,i)\big),\\
b'(i,j)&=&\dis b(j,i)+\ga(i,j),\\
c'(i,j)&=&\dis c(i,j)-\big(\ga(j,i)+\zeta(j,i)\big)
-\eta\big(\ga(i,j)+\zeta(i,j)\big)+\eta\big(b(i,j)-b(j,i)\big)\\
&=&\dis c(i,j)-\ga(j,i)-\eta\ga(i,j)+(1-\eta)\zeta(i,j)
+\eta\big(b(i,j)-b(j,i)\big)\\
&=&\dis c(i,j)-\ga(j,i)-\eta\ga(i,j)+\big(e(i,j)-e(j,i)\big)
+\eta\big(b(i,j)-b(j,i)\big),\\
d'(i,j)&=&\dis d(i,j)+\ga(i,j)+\big(e(i,j)-e(j,i)\big),\\
e'(i,j)&=&\dis e(j,i)-\ga(i,j).
\ec
An equivalent set of equations is
\bc\label{acrat3}
a'(i,j)&=&\dis a(i,j)+\eta\big(\ga(i,j)+\ga(j,i)\big),\\
b'(i,j)&=&\dis b(j,i)+\ga(i,j),\\
c'(i,j)+e'(i,j)+\eta b'(i,j)&=&\dis c(i,j)+e(i,j)+\eta b(i,j)
-\big(\ga(i,j)+\ga(j,i)\big),\\
d'(i,j)+e'(i,j)&=&\dis d(i,j)+e(i,j),\\
e'(i,j)&=&\dis e(j,i)-\ga(i,j).
\ec
Combining the first, third, and fourth equation, we see that
\be\ba{l}
(1-\eta)\ga'(i,j)=a'(i,j)+c'(j,i)-d'(i,j)+\eta b'(j,i)-e'(i,j)+e'(j,i)\\
\dis=a(i,j)+\eta\big(\ga(i,j)+\ga(j,i)\big)+c(j,i)+e(j,i)+\eta b(j,i)
-\big(\ga(j,i)+\ga(i,j)\big)-d(i,j)-e(i,j)\\
\dis=(1-\eta)\ga(i,j)-(1-\eta)\big(\ga(i,j)+\ga(j,i)\big)=-(1-\eta)\ga(j,i).
\ec
Thus, an equivalent set of equations for $a',\ldots,e'$ is
\bc\label{acrat4}
a'(i,j)+\eta\big(e'(i,j)+e'(j,i)\big)&=&\dis a(i,j)+\eta\big(e(i,j)
+e(j,i)\big),\\
b'(i,j)+e'(i,j)&=&\dis b(j,i)+e(j,i),\\
\ga'(i,j)&=&-\ga(j,i),\\
d'(i,j)+e'(i,j)&=&\dis d(i,j)+e(i,j),\\
e'(i,j)&=&\dis e(j,i)-\ga(i,j).
\ec
Combining the third and fifth equation yields (\ref{eqeq}).

What follows here is a check of (\ref{acrat}) with the formula for
$c'$ replaced by
\be
c'(i,j)=c(i,j)-(1-\eta)\oldga(j,i)-2\eta\oldga(i,j)+2\eta\zeta(i,j)
+\eta\big(b(i,j)-b(j,i)\big).
\ee
We got this wrong formula since our third formula in (\ref{ducofin2})
erroneously read
\be
c'(j,i)=-a'(i,j)+d'(i,j)-\eta b'(j,i)-\big(e(i,j)-e'(i,j)\big)
+\eta(e(j,i)-e'(i,j)).
\ee
To check our solution, we insert (\ref{acrat}) into (\ref{ducofin}).
This yields
\bc
a'(i,j)&=&\dis a(i,j)+\eta\big(\oldga(i,j)+\oldga(j,i)\big),\\
b'(i,j)&=&\dis b(j,i)+\oldga(i,j)-\zeta(i,j),\\
c'(i,j)&=&\dis c(i,j)-(1-\eta)\oldga(j,i)-2\eta\oldga(i,j)+2\eta\zeta(i,j)
+\eta\big(b(i,j)-b(j,i)\big),\\
d'(i,j)&=&\dis d(i,j)+\oldga(i,j)-\eta\zeta(i,j),\\
e'(i,j)&=&\dis e(j,i)-\oldga(i,j)+\zeta(i,j),\\
\ec
\bc
a'(j,i)&=&\dis a(j,i)+\eta\big(\oldga(j,i)+\oldga(i,j)\big),\\
b'(j,i)&=&\dis b(i,j)+\oldga(j,i)-\zeta(j,i),\\
c'(j,i)&=&\dis c(j,i)-(1-\eta)\oldga(i,j)-2\eta\oldga(j,i)+2\eta\zeta(j,i)
+\eta\big(b(j,i)-b(i,j)\big),\\
d'(j,i)&=&\dis d(j,i)+\oldga(j,i)-\eta\zeta(j,i),\\
e'(j,i)&=&\dis e(i,j)-\oldga(j,i)+\zeta(j,i),\\
\ec
\be\ba{l}\label{ducheck1}
\dis a(i,j)+\eta\big(\oldga(i,j)+\oldga(j,i)\big)+\eta\big(e(j,i)-\oldga(i,j)
+\zeta(i,j)+e(i,j)-\oldga(j,i)+\zeta(j,i)\big)\\
\dis\qquad=a(i,j)+\eta\big(e(i,j)+e(j,i)\big),\\[5pt]
\dis-(1-\eta)\big(e(j,i)-\oldga(i,j)+\zeta(i,j)\big)=a(i,j)-d(i,j)
-e(i,j)+\eta b(j,i)+c(j,i)+\eta e(j,i),\\[5pt]
\dis a(i,j)+\eta\big(\oldga(i,j)+\oldga(j,i)\big)-\big(d(i,j)+\oldga(i,j)
-\eta\zeta(i,j)\big)-\big(e(j,i)-\oldga(i,j)+\zeta(i,j)\big)\\
\dis\quad+\eta\big(b(i,j)+\oldga(j,i)-\zeta(j,i)\big)+c(j,i)
-(1-\eta)\oldga(i,j)-2\eta\oldga(j,i)+2\eta\zeta(j,i)\\
\dis\quad+\eta\big(b(j,i)-b(i,j)\big)+\eta\big(e(i,j)-\oldga(j,i)
+\zeta(j,i)\big)=-(1-\eta)e(i,j),\\[5pt]
\dis d(i,j)+\oldga(i,j)-\eta\zeta(i,j)+e(j,i)-\oldga(i,j)+\zeta(i,j)
=d(i,j)+e(i,j),\\[5pt]
\dis b(i,j)+\oldga(j,i)-\zeta(j,i)+e(i,j)-\oldga(j,i)+\zeta(j,i)=b(i,j)+e(i,j).
\ec
This can be simplified to
\be\ba{l}\label{ducheck2}
\dis \eta\big(\oldga(i,j)+\oldga(j,i)\big)+\eta\big(e(j,i)-\oldga(i,j)
+e(i,j)-\oldga(j,i)\big)\\
\dis\qquad=\eta\big(e(i,j)+e(j,i)\big),\\[5pt]
\dis-(1-\eta)\big(e(j,i)-\oldga(i,j)+\zeta(i,j)\big)=(1-\eta)\oldga(i,j)
-e(i,j)+\eta e(j,i),\\[5pt]
\dis a(i,j)+c(j,i)-d(i,j)+\eta b(j,i)+\eta\big(\oldga(i,j)+\oldga(j,i)\big)
-\big(\oldga(i,j)-\eta\zeta(i,j)\big)\\
\dis\quad-\big(e(j,i)-\oldga(i,j)+\zeta(i,j)\big)+\eta\big(\oldga(j,i)
-\zeta(j,i)\big)-(1-\eta)\oldga(i,j)-2\eta\oldga(j,i)+2\eta\zeta(j,i)\\
\dis\quad+\eta\big(e(i,j)-\oldga(j,i)+\zeta(j,i)\big)=-(1-\eta)e(i,j),\\[5pt]
\dis -\eta\zeta(i,j)+\zeta(i,j)=e(i,j)-e(j,i),\\[5pt]
\dis 0=0.
\ec
and further to 
\be\ba{l}\label{ducheck3}
\dis 0=0,\\[5pt]
\dis-(1-\eta)e(j,i)-\big(e(i,j)-e(j,i)\big)=-e(i,j)+\eta e(j,i),\\[5pt]
\dis (1-\eta)\oldga(i,j)-(1-2\eta)\oldga(i,j)-\eta\oldga(j,i)
-(1-\eta)\zeta(i,j)-e(j,i)+\eta e(i,j)=-(1-\eta)e(i,j),\\[5pt]
\dis -\eta\zeta(i,j)+\zeta(i,j)=(1-\eta)\zeta(i,j),\\[5pt]
\dis 0=0,
\ec
and
\be\ba{l}\label{ducheck4}
\dis 0=0,\\[5pt]
\dis 0=0,\\[5pt]
\dis \eta\oldga(i,j)-\eta\oldga(j,i)-(1-\eta)\big(e(i,j)-e(j,i)\big)-e(j,i)
=-e(i,j),\\[5pt]
\dis 0=0,\\[5pt]
\dis 0=0.
\ec
Dividing out a factor $\eta$, we end up with
\be
\oldga(i,j)-\oldga(j,i)+e(i,j)-e(j,i)=0.
\ee
This is false, since $e$ might be symmetric while $\oldga$ is not.

We must check our equations for $c'$ again. We start with the third
equation in (\ref{ducofin}):
\be
a'(i,j)-d'(i,j)-e'(i,j)+\eta b'(j,i)+c'(j,i)+\eta e'(j,i)=-(1-\eta)e(i,j).
\ee
This can be rewritten as
\be
c'(j,i)=-a'(i,j)+d'(i,j)-\eta b'(j,i)+e'(i,j)-\eta e'(j,i)-(1-\eta)e(i,j),
\ee
and
\be\label{cac}
c'(j,i)=-a'(i,j)+d'(i,j)-\eta b'(j,i)-\big(e(i,j)-e'(i,j)\big)
+\eta\big(e(i,j)-e'(j,i)\big).
\ee
If we believe in (\ref{edif}) and (\ref{ducoabd}), this implies
\bc
c'(j,i)&=&\dis-\big(a(i,j)+\eta\big(\oldga(i,j)+\oldga(j,i)\big)\big)
+\big(d(i,j)+\oldga(i,j)-\eta\zeta(i,j)\big)\\[5pt]
&&\dis-\eta\big(b(i,j)+\oldga(j,i)-\zeta(j,i)\big)-\big(\oldga(i,j)
-\eta\zeta(i,j)\big)+\eta\big(\oldga(j,i)-\zeta(j,i)\big)\\[5pt]
&=&\dis c(j,i)-\big(a(i,j)+c(j,i)-d(i,j)+\eta b(j,i)\big)
-\eta\big(b(i,j)-b(j,i)\big)\\
&&\dis-\eta\big(\oldga(i,j)+\oldga(j,i)\big)\\[5pt]
&=&\dis c(j,i)-(1-\eta)\oldga(i,j)-\eta\big(b(i,j)-b(j,i)\big)\\
&&\dis-\eta\big(\oldga(i,j)+\oldga(j,i)\big)\\[5pt]
&=&\dis c(j,i)-\oldga(i,j)-\eta\oldga(j,i)-\eta\big(b(i,j)-b(j,i)\big).
\ec
let us compare (\ref{cac}) with the third formula in (\ref{ducofin2}).
The latter reads
\be\label{thirdu}
c'(j,i)=-a'(i,j)+d'(i,j)-\eta b'(j,i)-\big(e(i,j)-e'(i,j)\big)
+\eta(e(j,i)-e'(i,j)).
\ee
We spot a difference in the last term. The third formula in (\ref{ducofin2})
is wrong here!}

\noi
{\bf Remark} The original motivation for the calculations in this
appendix was to generalize the self-duality in
Proposition~\ref{P:cvdu}~(a) to the case that $q(i,j)\neq
q(j,i)$. Surprisingly, this does not work, unless $r=0$. In
particular, if $m=0$ but $r,s>0$, one can check from (\ref{newrates})
that a $(q,r,s,0)$-CVP can only be dual to another $(q',r',s',0)$-CVP
if $q=q'$ and $q(i,j)=q(j,i)$.

\detail{Indeed, the formula for $e'$ forces $\ga=0$. Since $m=0$ one
has $c=0$. Since $\ga=0$ one has $\eta\neq 0$. Now the requirement
$c'=0$ forces $q(i,j)-q(j,i)=0$.}

\section{Local mean field limits}

\subsection{Contact voter processes}\label{A:cv}

{\bf Idea of proof of Claim~\ref{C:cv}} The process $\ov X^{(N)}$
jumps as \be\ba{lcl} \dis x\to x+\ffrac{1}{N}\de_i&\quad\mbox{with
rate}\quad&\dis\nu^{(N)}N^{-1}(r^{(N)}+s)\sum_{j:\,j\neq
i}q(j,i)x(j)N(N-x(i)N)\\
&&\dis+(1-\nu^{(N)})(N-1)^{-1}(r^{(N)}+s)x(i)N(N-x(i)N)\\[10pt] \dis
x\to x-\ffrac{1}{N}\de_i&\quad\mbox{with
rate}\quad&\dis\nu^{(N)}N^{-1}r^{(N)}\sum_{j:\,j\neq
i}q(j,i)(N-x(j)N)x(i)N\\
&&\dis+(1-\nu^{(N)})(N-1)^{-1}r^{(N)}(N-x(i)N)x(i)N\\[5pt]
&&\dis+mx(i)N.  \ec It follows that the process started in $\ov
X^{(N)}_0=x$ satisfies \bc \dis\lim_{t\to 0}t^{-1}E[\ov
X^{(N)}_t(i)-x(i)]&=&\dis\nu^{(N)}r^{(N)}\sum_{j:\,j\neq
i}q(j,i)(x(j)-x(i))\\ &&\dis+\nu^{(N)}s\sum_{j:\,j\neq
i}q(j,i)x(j)(1-x(i))\\
&&\dis+(1-\nu^{(N)})s(1-N^{-1})^{-1}x(i)(1-x(i))\\[5pt]
&&\dis-mx(i)\\[10pt] \dis\lim_{t\to 0}t^{-1}E[(\ov
X^{(N)}_t(i)-x(i))^2]&=&\dis\nu^{(N)}(2r^{(N)}+s)N^{-1}\sum_{j:\,j\neq
i}q(j,i)x(j)(1-x(i))\\
&&\dis+(1-\nu^{(N)})(2r^{(N)}+s)(N-1)^{-1}x(i)(1-x(i))\\[5pt]
&&\dis+mx(i)N^{-1}.  \ec Under the assumptions (\ref{raco}), this
simplifies to \bc \dis\lim_{t\to 0}t^{-1}E[\ov
X^{(N)}_t(i)-x(i)]&=&\dis\sum_{j:\,j\neq i}q(j,i)(x(j)-x(i))\\
&&\dis+sx(i)(1-x(i))-mx(i)+O(N^{-1}),\\[5pt] \dis\lim_{t\to
0}t^{-1}E[(\ov X^{(N)}_t(i)-x(i))^2]&=&\dis 2r x(i)(1-x(i))+O(N^{-1})
\ec as $N\to\infty$. This, plus an estimate on a higher moment, can be
used to show that the limiting process $\Xc$ started in $\Xc_0=x$
satisfies \bc\label{todif} \dis\lim_{t\to
0}t^{-1}\big(E[f(\Xc_t)]-f(x)\big)&=&\dis\sum_{i\neq
j}q(j,i)(x(j)-x(i))\dif{x(i)}f(x)+s\sum_ix(i)(1-x(i))\dif{x(i)}f(x)\\[5pt]
&&\dis-m\sum_ix(i)\dif{x(i)}f(x)+r\sum_ix(i)(1-x(i))\diff{x(i)}f(x)
\ec for all twice continuously differentiable $f$ depending on
finitely many coordinates. In the right-hand side of (\ref{todif}) we
reckognize the generator of the diffusion process in
(\ref{SSMsde}). More formally, the calculations above can be used to
show that the processes $\ov X^{(N)}$ are tight, and each weak limit
point satisfies the martingale problem for the operator in
(\ref{todif}). Standard theory then yields the convergence in
(\ref{toSSM}).

To prove the convergence in (\ref{toBPS}), set
\bc
\rho^{(N)}&:=&r^{(N)}+\frac{\eps}{1+\eps}s,\\
\bet&:=&\frac{1}{1+\eps}s,\\
\de&:=&m-\frac{\eps}{1+\eps}s,
\ec
and observe that the assumptions (\ref{raco}) imply that
\bc\label{ram}
\rho^{(N)}\nu^{(N)}&\to&1,\\
\rho^{(N)}N^{-1}&\to&r.
\ec
The process $\ov Y^{(N)}$ jumps as
\be\ba{lcl}
\dis y\to y-\de_i+\de_j&\quad\mbox{with rate}\quad&
\dis\nu^{(N)}N^{-1}\rho^{(N)}q(i,j)y(i)(N-y(j)),\\[10pt]
\dis y\to y-\de_i-\de_j&\quad\mbox{with rate}\quad&
\dis\eps\nu^{(N)}N^{-1}\rho^{(N)}q(i,j)y(i)y(j),\\[10pt]
\dis y\to y+\de_i&\quad\mbox{with rate}\quad&
\dis\nu^{(N)}N^{-1}\bet\sum_{j:\,j\neq i}q(j,i)y(j)(N-y(i))\\
&&\dis+(1-\nu^{(N)})(N-1)^{-1}\bet y(i)(N-y(i)),\\[10pt]
\dis y\to y-\de_i&\quad\mbox{with rate}\quad&
\dis(1-\eps)\nu^{(N)}N^{-1}\rho^{(N)}\sum_{j:\,j\neq i}q(i,j)y(i)y(j)\\
&&\dis+(1-\eps)(1-\nu^{(N)})(N-1)^{-1}\rho^{(N)}y(i)(y(i)-1)\\[5pt]
&&\dis+\eps\nu^{(N)}N^{-1}\bet\sum_{j:\,j\neq i}q(j,i)y(j)y(i)\\
&&\dis+\eps(1-\nu^{(N)})(N-1)^{-1}\bet y(i)(y(i)-1)\\[5pt]
&&\dis+\de y(i),\\[10pt]
\dis y\to y-2\de_i&\quad\mbox{with rate}\quad&
\dis\eps(1-\nu^{(N)})(N-1)^{-1}\rho^{(N)}y(i)(y(i)-1).
\ec
In the limit $N\to\infty$, this simplifies to
\be\ba{lcl}
\dis y\to y-\de_i+\de_j&\quad\mbox{with rate}\quad&\dis q(i,j)y(i),\\[10pt]
\dis y\to y+\de_i&\quad\mbox{with rate}\quad&
\dis\ffrac{1}{1+\eps}sy(i),\\[10pt]
\dis y\to y-\de_i&\quad\mbox{with rate}\quad&
\dis(1-\eps)ry(i)(y(i)-1)+(m-\ffrac{\eps}{1+\eps}s)y(i),\\[10pt]
\dis y\to y-2\de_i&\quad\mbox{with rate}\quad&\dis\eps ry(i)(y(i)-1),
\ec
plus terms of order $N^{-1}$.\qed

\subsection{Super random walks}\label{A:srw}

{\bf Idea of proof of Claim~\ref{C:cv}} The proof of (\ref{XtoSRW}) is
very similar to the proof of (\ref{toSSM}). Indeed, the process
$\ov X^{(N)}$ jumps as
\be\ba{lcl}
\dis x\to x+\ffrac{\phi}{\sqrt N}\de_i&\quad\mbox{with rate}\quad&
\dis\nu^{(N)}N^{-1}s^{(N)}\sum_{j:\,j\neq i}q(j,i)x(j)\phi^{-1}
\sqrt N(N-x(i)\phi^{-1}\sqrt N)\\
&&\dis+(1-\nu^{(N)})(N-1)^{-1}s^{(N)}x(i)\phi^{-1}\sqrt N(N-x(i)
\phi^{-1}\sqrt N)\\[10pt]
\dis x\to x-\ffrac{\phi}{\sqrt N}\de_i&\quad\mbox{with rate}\quad&
\dis m^{(N)}x(i)\phi^{-1}\sqrt N.
\ec
It follows that the process started in $\ov X^{(N)}_0=x$ satisfies
\bc
\dis\lim_{t\to 0}t^{-1}E[\ov X^{(N)}_t(i)-x(i)]&=&
\dis\nu^{(N)}s^{(N)}\sum_{j:\,j\neq i}q(j,i)x(j)(1-x(i)\phi^{-1}N^{-1/2})\\
&&\dis+(1-\nu^{(N)})s^{(N)}(x(i)-\phi^{-1}N^{-1/2}x(i)^2)(1-N^{-1})^{-1}\\[5pt]
&&\dis- m^{(N)}x(i)\\[10pt]
\dis\lim_{t\to 0}t^{-1}E[(\ov X^{(N)}_t(i)-x(i))^2]&=&
\dis\phi N^{-1/2}\nu^{(N)}s^{(N)}\sum_{j:\,j\neq i}q(j,i)x(j)(1-x(i)
\phi^{-1}N^{-1/2})\\
&&\dis+\phi N^{-1/2}(1-\nu^{(N)})s^{(N)}(x(i)
-\phi^{-1}N^{-1/2}x(i)^2)(1-N^{-1})^{-1}\\[5pt]
&&\dis+\phi N^{-1/2}m^{(N)}x(i).
\ec
By (\ref{raco2}), for large $N$ this simplifies to
\bc
\dis\lim_{t\to 0}t^{-1}E[\ov X^{(N)}_t(i)-x(i)]&=&
\dis\sum_{j:\,j\neq i}q(j,i)x(j)+(\bet-1)x(i)-\ga x(i)^2\\[5pt]
&=&\dis\sum_{j:\,j\neq i}q(j,i)(x(j)-x(i))+\bet x(i)-\ga x(i)^2,\\[5pt]
\dis\lim_{t\to 0}t^{-1}E[(\ov X^{(N)}_t(i)-x(i))^2]&=&\dis\al x(i),
\ec
plus terms of order $N^{-1/2}$. (Here we have used that
$\sum_{j:\;j\neq i}q(j,i)=1$.) It follows that $\ov X^{(N)}$
converges to the process with generator
\bc
Gf(z)&=&\dis\sum_{i\neq j}q(j,i)(z(j)-z(i))\dif{z(i)}
+\ffrac{1}{2}\al\sum_iz(i)\diff{z(i)}\\
&&\dis+\bet\sum_iz(i)\dif{z(i)}-\ga\sum_iz(i)^2\dif{z(i)},
\ec
i.e., the process in (\ref{SRWsde}).\qed


\begin{thebibliography}{BCGH95}

\bibitem[AS05]{AS05}
S.R.~Athreya and J.M.~Swart.
Branching-coalescing particle systems.
{\em Prob.\ Theory Relat.\ Fields}~131(3), 376--414, 2005.

\bibitem[EK86]{EK}
S.N.~Ethier and T.G.~Kurtz.
{\em Markov Processes; Characterization and Convergence.}
John Wiley \& Sons, New York, 1986.

\bibitem[HWS05]{HW05}
M.~Hutzenthaler and A.~Wakolbinger.
Ergodic behaviour of locally regulated branching populations.
ArXiv math.PR/0509612.

\bibitem[Lig85]{Lig85}
T.M.~Liggett.
{\em Interacting Particle Systems}.
Springer, New York, 1985.

\bibitem[SL95]{SL95}
A.~Sudbury and P.~Lloyd.
Quantum operators in classical probability theory. II: The concept of
duality in interacting particle systems.
{\em Ann.\ Probab.}~23(4), 1816--1830, 1995.

\bibitem[SL97]{SL97}
A.~Sudbury and P.~Lloyd.
Quantum operators in classical probability theory. IV: Quasi-duality and
thinnings of interacting particle systems
{\em Ann.\ Probab.}~25(1), 96--114, 1997. 

\bibitem[SU86]{SU86}
T.~Shiga and K.~Uchiyama.
Stationary states and their stability of the stepping stone model
involving mutation and selection.
{\em Probab.\ Theory Relat.\ Fields} 73, 87--117, 1986.

\bibitem[Sud00]{Sud00}
A.~Sudbury.
Dual families of interacting particle systems on graphs.
{\em J.\ Theor.\ Probab.}~13(3), 695--716, 2000.

\end{thebibliography}
\end{document}